




\documentstyle{amsppt}
\define\z{{\bold Z }}
\define\cc{{\bold C}}
\define\r{{\bold R}}
\define\f{{\bold F}}

\define\Het{H_{et}}
\define\Ker{\text{Ker}\ }
\define\im{\text{Im}\ }
\define\Pic{\text{Pic}\ }
\magnification1200

\topmatter

\title
ON THE BRAUER GROUP OF REAL AlGEBRAIC SURFACES
\endtitle

\author
Viacheslav V. Nikulin
\endauthor

\address
Steklov Mathematical Institute,
ul. Vavilova 42, Moscow 117966, GSP-1, Russia.
\endaddress

\email
slava\@nikulin.mian.su
\endemail

\dedicatory
Dedicated to Professor Igor R.
Shafarevich on the occasion of his seventieth birthday
\enddedicatory

\abstract
     Let $X$ be a real projective algebraic manifold, $s$ numerates
connected components of $X(\r)$ and $_2Br(X)$ the subgroup of
elements of order $2$ of cohomological Brauer group $Br(X)$.

   We study the natural homomorphism
$\xi : _2Br(X) \to (\z/2)^s$ and prove that $\xi$ is epimorphic if
$H^3(X(\cc)/G;\z/2)\to H^3(X(\r);\z/2)$
is injective. Here $G=Gal(\cc/\r)$.

   For an algebraic surface $X$ with $H^3(X(\cc)/G;\z /2)=0$ and
$X(\r)\not=\emptyset$, we
give a formula for
$\dim_2Br(X)$.

    As a corollary, for a real Enriques surfaces $Y$, the
$\xi$ is epimorphic and $\dim_2Br(Y)=2s-1$ if both liftings of
the antiholomorphic involution of $Y$ to the universal covering $K3$-
surface  $X$ have non-empty sets of real points (this is the general case).
For this case, we also give a formula for the number $s_{nor}$ of
non-orientable components of $Y$ which is very important for
the topological classification of real Enriques surfaces.
\endabstract

\endtopmatter

\document

\head
\S 0. Formulation of basic results
\endhead

In the paper of R. Sujatha and the author
\cite{N-S}, the Brauer group of a real Enriques surface was studied.
Here we continue the study of Brauer group with
the remark that most of the  results of these paper generally valid for
an arbitrary smooth projective real algebraic surface.

Let $X$ be a projective algebraic
variety over the field $\r$ of real numbers.
Let
$$
Br^\prime (X)=H^2_{et}(X; {\bold G}_m)
$$
denote the cohomological Brauer group
of $X$. See a definition in the book of J. Milne \cite{Mi}, for example.
We mention that the cohomological Brauer group is very closely related with
the more interesting
classical Brauer group $Br(X)$ classifying Azumaya algebras over $X$
(see the papers of A. Grothendieck \cite{Gr2}
and the book \cite{Mi}). For example,
it is known that $Br(X)\subset Br^\prime (X)$. For  curves and smooth
surfaces it gives an isomorphism. But we will only consider  here the
cohomological Brauer group $Br^\prime (X)$.

Let $X(\r )$ denote the space of
$\r$-rational points of $X$ with the Euclidean topology and $s$ denote the
number of
real connected components of this space.
Let $_2Br^\prime (X)$ denote the group of
elements of order two in $Br^\prime (X).$
If $P\in X(\r )$ is
a real point of $X$, we get a natural map
$_2Br^\prime (X)\to _2Br^\prime (P)\cong \z /2$.
It is shown in \cite{CT-P} that this map
does not depend from a choice of the
point $P$ in  a connected component of $X(\r )$. Thus, the canonical
 map
$$
_2Br^\prime (X)\to (\z /2)^s
\tag0--1
$$
is defined.

We mention that a studying of the map (0--1)
and a description of the Brauer
group of $X$ is very important for calculation of a such interesting
group connected  with $X$ as the Witt group $W(X)$.
See R. Sujatha \cite{Su}.

It is shown in the paper of J.-L. Colliot-Th\'el\`ene and R. Parimala
\cite{CT-P}
that
the map (0--1) is epimorphic if $X$ is
a smooth projective algebraic surface and
$H^3(X(\cc ); \f_2)=0$ (here $\f _2=\z /2$). This is a generalization of the
old result of E. Witt \cite{W} about curves
(compare with Remark 1.8 below).
There doesn't seem to be any known example of a surface where the map  (0--1)
fails to be epimorphic.

This paper is devoted to studying of this map (0--1) and also calculating of
$\dim~_2Br^\prime (X)$.
Our idea is to interpret $_2Br^\prime (X)$ and the
map (0--1) purely topologically, and apply to them topological
considerations.

We prove here the following basic result where $G$ is the group
of order two
generated by the antiholomorphic involution $g$ of $X(\cc)$ defined by the
structure of real algebraic variety on $X$.

\proclaim{Theorem 0.1} Let $X/\r $ be an algebraic projective manifold
(smooth) over the field $\r$ of real numbers.

Then, the
homomorphism  (0--1) is  epimorphic if   $H^3(X(\cc )/G; \f_2)=0$.
More generally, the homomorphism (0--1) is epimorphic if the kernel of the
homomorphism
$$
i^\ast : H^3(X(\cc)/G;\f_2)\to H^3(X(\r);\f_2)
$$
is equal to zero. Here $i :X(\r)\subset X(\cc)/G$ denote the embedding.
\endproclaim

\demo{Proof} See Theorem 1.6 below. In fact, in
Theorem 1.6,  we give a precise
topological  obstruction to epimorphicity of the map (0--1).
This obstruction is zero if the kernel of the homomorphism $i^\ast$ above is
zero.

We mention that for smooth curves $X$ the group
$H^3(X(\cc )/G;\f_2)=0$, thus the
map (0--1) is  epimorphic (it is well-known, compare with \cite{W}).
For surfaces $X$
the group $H^3(X(\r);\f_2)=0$, and  $\Ker i^\ast =H^3(X(\cc )/G;\f_2)$.
\enddemo

Now, let us show that, from Theorem 0.1, the
result of J.-L. Colliot-Th\'el\`ene and R. Parimala mentioned above follows.

Let  $X$ be a smooth projective algebraic
surface and  $H^3(X(\cc ); \f _2)=0$.
Then, by  Poincar\'e duality, we have
$H_1(X(\cc);\f_2)=0$. If $X(\r)\not=\emptyset$,
then  any loop in $X(\cc )/G$ with the beginning on $X(\r)$ has a lifting
to a loop  on $X(\cc)$.
It follows that  the canonical homomorphism $H_1(X(\cc);\f_2)
\to  H_1(X(\cc)/G; \f_2)$ is epimorphic. Thus, $H_1(X(\cc)/G;\f_2)=0$.
For the dimension $2$ the quotient space
$X(\cc)/G$ is homeomorphic to a smooth compact $4$-dimensional
manifold.   By  Poincar\'e duality, we then get that
$H^3(X(\cc)/G; \f_2)=0$. This proves the  statement.

\vskip 5pt

We apply the Theorem 0.1 to real Enriques surfaces.

By a complex Enriques surface $Y$ over $\cc$ ,
we mean a non-singular minimal projective algebraic
surface $Y/\cc$ such that the invariants
$\kappa(Y)=p_g(Y)=q(Y)=0$.
These are equivalent to irregularity $q(Y)=0$ and $2K_Y=0$ but
$K_Y\not=0$ where $K_Y$  is the canonical class of $Y$.
One may find all information about
Enriques surfaces we need in the books \cite{A} and
\cite{C-D}.

By a real Enriques surface $Y/\r $, we mean a projective
algebraic surface $Y/\r $ such that $Y\otimes _\r \cc $ is a complex Enriques
surface. Universal covering complex surface of an Enriques surface
$Y(\cc)$ is a $K3$-surface $X(\cc)$ (see \cite{A} and \cite{C-D})
which twice covers the Enriques surface
$Y(\cc)$. We denote by $\tau$ the holomorphic involution on $X(\cc)$ of
this covering. There are precisely two liftings $\sigma$ and
$\tau\sigma $ on $X(\cc)$ of the
antiholomorphic involution $\theta$ of
$Y(\cc)$ corresponding to the real structure on $Y$. Besides, one can see
very easily that if
$Y(\r)\not=\emptyset$, then both $\sigma$ and
$\tau\sigma$ are antiholomorphic involutions of $X(\cc)$. Thus,
$\sigma$ and $\tau\sigma$ define two real structures
$X_\sigma$ and $X_{\tau\sigma}$ on the $K3$-surface $X$. We
denote by
$$
X_\sigma(\r)=X(\cc)^\sigma ,\ \ \ \
X_{\tau\sigma}(\r)=X(\cc)^{\tau\sigma}
$$
the real parts of the real $K3$-surfaces $X_\sigma$ and
$X_{\tau\sigma}$ corresponding
to these real structures respectively.
Since $\tau$ has no fixed points on $X(\cc)$, it follows that the
sets  $X_\sigma(\r)$ and $X_{\tau\sigma}(\r)$ have an empty intersection.
 From the Theorem 0.1, we get

\proclaim{Corollary 0.2}
Let $Y$ be a real Enriques surface with the antiholomorphic involution
$\theta$, and the real part $Y(\r)\not=\emptyset$.
Suppose that the real parts $X_\sigma(\r)$ and
$X_{\tau\sigma}(\r)$ of both liftings $\sigma$ and $\tau\sigma$ of
$\theta$ to the universal covering $K3$-surface $X(\cc)$ are non-empty.

Then the canonical map (0--1) corresponding to the real Enriques surface
$Y$ is epimorphic.
\endproclaim

\demo{Proof}
$$
Y(\cc )/\{ \text{id} _{Y(\cc )},\theta \} =
X(\cc )/\{ \text{id} _{X(\cc )}, \tau , \sigma ,\tau \sigma \} =
$$
$$
(X(\cc )/\{ \text{id}_{X(\cc )}, \sigma \} )/\{ \text{id} , \tau \sigma
\mod \{ \text{id}_{X(\cc )}, \sigma \} \} .
$$
Here
involutions $\sigma $ and
$\tau \sigma ~\mod \{ \text{id}_{X(\cc )}, \sigma \} $
have non-empty sets of fixed
points because real parts of both involutions $\sigma$ and $\tau \sigma $ of
$X(\cc )$ are non-empty and are not coincided.  Since for  a $K3$-surface $X$,
the group
$H_1(X(\cc );\f _2)=0$, it follows like above that
$H_1(X(\cc )/\{ \text{id}_{X(\cc )}, \sigma \};\f_2)=0$ and
$H_1(Y(\cc )/\{ \text{id} _{Y(\cc)};\theta \};\f _2)=0$.
The  topological space
$Y(\cc )/\{ \text{id} _{Y(\cc )};\theta \}$ is isomorphic to a
smooth compact  $4$-dimensional manifold.
By Poincar\'e
duality, then
$H^3(Y(\cc )/\{ \text{id} _{Y(\cc )}, \theta \};\f _2)=0$.
By Theorem 0.1,
the map  (0--1) is epimorphic for the real Enriques surface $Y$.
\enddemo

The same considerations show that if one of involutions $\sigma $ or
$\tau \sigma $ has an empty set of real points, then
$$
H_1(Y(\cc )/\{ \text{id} _{Y(\cc )},\theta \};\f _2)=\f _2
\ \text{and} \
H^3(Y(\cc )/\{ \text{id} _{Y(\cc )},\theta \};\f _2)=\f _2.
$$
-- If $X_\sigma (\r )\not= \emptyset$ but
$X_{\tau\sigma}(\r)=\emptyset$
then  the surface
$
X(\cc )/\{ \text{id} _{X(\cc )},  \sigma \}
$
is a 2-sheeted universal covering of the surface
$Y(\cc )/\{ \text{id} _{Y(\cc )},\theta \}$.
Thus,  the Corollary 0.2 gives precisely the case when the Theorem 0.1
may be applied to real Enriques surfaces.

We discuss in Remark 1.7 below a chance of constructing a counterexample to
epimorphisity of the map (0--1) using real Enriques surfaces $Y$ above with
$H^3(Y(\cc )/G;\f_2)=\f_2$ (equivalently, with
$X_\sigma(\r)\not=\emptyset$ but $X_{\tau\sigma}(\r)=\emptyset$). It is
not difficult to show that real Enriques surfaces with the condition
$H^3(Y(\cc )/G;\f_2)=\f_2$ do
exist. For the most part of real Enriques surfaces (from the point of view
of
the number of connected components of the moduli space) both involutions
$\sigma$ and $\tau\sigma$ have a non-empty set of real points. But for
some real Enriques surfaces one of these involutions may have an empty set
of real points.

\vskip 10pt

The following results are
devoted to a calculation of the dimension of \'etale
cohomology  groups with coefficients $\f_2$, and $_2Br^\prime (X)$.

We begin with the following general remark about a real algebraic
variety $X$.

We recall that the Kummer sequence
$$
0 \to \mu_2 \to {\bold  G}_m \to  {\bold  G}_m \to 0
\tag0--2
$$
yields the exact sequence
$$
0 \to \Pic ~X/2\Pic ~X \to H^2_{et}(X; \mu _2) \to _2~Br^\prime (X)
\to 0.
\tag0--3
$$

If $X(\r )$ is non-empty, $\Pic X=(\Pic (X \otimes \cc))^G$
(it is well-known \cite{Ma} and not difficult to see).
Thus, from (0--3), we have
$$
\dim~_2Br^\prime (X)=\dim~H^2_{et}(X;\mu_2)-
\dim~(\Pic (X \otimes \cc ))^G/2(\Pic (X \otimes \cc ))^G.
\tag0--4
$$
The dimension of the \'etale cohomology group $H^2_{et}(X;\mu_2)=
H^2_{et}(X;\f _2)$ is
estimated using the Serre-Hochschild spectral sequence where
$G=Gal(\cc /\r )$,
$$
E_2^{p,q}=H^p(G;H^q_{et}(X\otimes \cc ;\f _2))
\Longrightarrow H^{p+q}_{et}(X;\f _2),
\tag0--5
$$
where for a complex manifold $X\otimes \cc$ we have
$H^q_{et}(X\otimes \cc ;\f_2)=H^q(X(\cc);\f_2)$
(see \cite{Mi}, for example).

In the \S 2, we prove the following results which
permit to calculate the dimension of the \'etale cohomology
groups with coefficients  $\f_2$ and the 2-torsion of the
Brauer group for surfaces satisfying to the
condition of Theorem 0.1. These  results show that the class of real smooth
projective surfaces $X$ satisfying to  the condition of
Theorem 0.1 is very nice (easy to work with).

\proclaim{Theorem 0.3} Let $X/\r$ be a real smooth projective algebraic
surface such that
$X(\r )\not=\emptyset$
and $H^3(X(\cc )/G;\f_2)=0$. Then the Serre-Hochschild spectral
sequence  (0--5) degenerates and
$$
\dim~H^0_{et}(X;\f _2)=1;
$$
$$
\dim~H^1_{et}(X;\f _2)=\dim~H^1(X(\cc );\f _2)+1;
$$
$$
\dim~H^2_{et}(X;\f _2)=\dim~H^2(X(\cc );\f _2)^G+\dim~H^1(X(\cc );\f _2)+1;
$$
$$
\split
\dim~H^3_{et}(X;\f _2)=2\dim~H^2(X(\cc );\f _2)^G & -\dim~H^2(X(\cc );\f _2)\\
&+2\dim~H^1(X(\cc );\f _2)+1
\endsplit
$$
$$
\split
\dim~H^k_{et}(X;\f _2)=2\dim~H^2(X(\cc );\f _2)^G &-\dim~H^2(X(\cc );\f _2)\\
&+2\dim~H^1(X(\cc );\f _2)+2
\endsplit
$$
for $k\ge 4$.
\endproclaim

Using Theorem 0.3 and (0--4), (0--5), we get

\proclaim{Theorem 0.4} Let $X/\r $ be a real smooth
projective algebraic surface such that $X(\r )\not= \emptyset $
and $H^3(X(\cc )/G;\f _2)=0$.

Then
$$
\dim~_2Br^\prime (X)=2s-1+h^{2,0}(X(\cc ))+h^{1,1}_-(X(\cc ))
- \rho _+(X\otimes \cc ).
$$
Here $h^{1,1}_-(X(\cc ))=\dim~H^{1,1}_-(X(\cc ))$ where
$$
H^{1,1}_-(X(\cc ))=\{ x\in H^{1,1}(X(\cc )) \mid g(x)=-x \}
$$
is the set of potentially real algebraic cycles. And
$\rho_+ (X\otimes \cc )=\dim~(\Pic (X \otimes \cc )\otimes \cc )^G$.
The characteristic class map gives
an injection of $(\Pic (X\otimes \cc )\otimes \cc)^G$ to
$H^{1,1}_-(X(\cc))$.
\endproclaim

For an Enriques surface, $h^{2,0}(Y(\cc))=0$ and all cycles are algebraic.
For a real Enriques surface $Y$, we have seen above that  the
condition $H^3(Y(\cc )/G;\f _2)=0$ is equivalent to
the condition of Corollary 0.2. Thus, we get

\proclaim{Corollary 0.5}
Let $Y$ be a real Enriques surface with the
antiholomorphic involution $\theta $ and
the real part $Y(\r)\not= \emptyset$.
Suppose that  the real parts $X_\sigma(\r)$ and
$X_{\tau\sigma}(\r)$ of both liftings $\sigma$ and $\tau\sigma$ of
$\theta$ to the universal covering $K3$-surface $X(\cc)$ are not empty.

Then the Serre-Hochschild spectral sequence (0--5) degenerates and
$$
\dim~_2Br^\prime (Y)=2s-1
$$
where $s$ is the number of real connected components of $Y(\r )$.
\endproclaim

We mention that for a real rational surface $Z$ with a non-empty set of
real points $Z(\r)$ the same results were known:
The map (0--1) is epimorphic
and $\dim~_2Br^\prime (Z)=2s-1$. It is also known that the Witt group
$W(Z)\cong (\z) ^s\oplus (\z/2)^{s-1}$. See \cite{Su}.
Perhaps,the last result about the Witt group
also valid for real Enriques surfaces with non-empty sets
$X_\sigma (\r)$ and $X_{\tau\sigma}(\r)$.

Using results of \cite{N-S}, we may prove some additional results
about Brauer groups of real Enriques surfaces which also valid if
one of the sets $X_\sigma (\r), X_{\tau\sigma}(\r)$ is empty.

In \cite{N-S}, the important invariants $b(Y)$ and $\epsilon(Y)$ of
a real Enriques surface $Y$ with an antiholomorphic involution $\theta$
were introduced. Here
$$
b(Y)=
\dim~H^2(Y(\cc);\f_2)^\theta -
\dim~(\Pic~Y\otimes \cc)^\theta /2(\Pic~Y\otimes \cc)^\theta + 1.
$$
The invariant $\epsilon(Y)=1$  if the differential $d_2^{0,2}$ of
the Hochschild-Serre spectral sequence (0--5) is zero, and
$\epsilon(Y)=0$ otherwise. We have the following results from
\cite{N-S} about these invariants:
$$
\dim~_2Br^\prime (Y)=b(Y)+\epsilon(Y)\ \ \text{if}\ \ Y(\r)\not=\emptyset,
\tag0--6
$$
$$
b(Y)\ge 2s-2 \ \ \text{for any $Y$}.
\tag0--7
$$
Thus, by (0--6) and (0--7), for any real Enriques surface $Y$,
$$
\dim~_2Br^\prime (Y)\ge 2s-2+\epsilon(Y).
\tag0--8
$$
Additionally to Corollary 0.5 and (0--6)---(0--8), we prove

\proclaim{Theorem 0.6} Let $Y$ be a real Enriques surface.
Then:

(i) The inequality (0--7) is an equality, i.e.
$b(Y)=2s-2$,
iff the Hochschild--Serre spectral sequence (0--5) degenerates. In
particular, by Corollary 0.5, $b(Y)=2s-2$ if $X_\sigma (\r)\not=
\emptyset$ and $X_{\tau\sigma}(\r)\not=\emptyset$.

(ii) $\dim~_2Br^\prime (Y)=2s-1$ if the Hochschild--Serre spectral sequence
(0--5) degenerates and $Y(\r)\not=\emptyset$. In particular, it is true if
$X_\sigma (\r)\not=\emptyset $ and $X_{\tau\sigma}(\r)\not=\emptyset$.

(iii) $\dim~_2Br^\prime (Y)\ge 2s-1$.
\endproclaim

In \S 3, we give an application of results of \cite{N-S} and
Corollary 0.5 and Theorem 0.6
to a topological studying of real Enriques surfaces $Y$.
Let $Y$ be a real Enriques surface with the antiholomorphic
involution $\theta$ and $Y(\r)\not=\emptyset$. Let $s_{nor}$ be
the number of non-orientable connected components of $Y(\r)$.
We denote by
$$
\Gamma =\{ \text{id}, \tau, \sigma, \tau\sigma \}\cong (\z/2)^2
$$
the group acting on the $K3$-surface $X(\cc)$ (we use notation
above). Let us suppose that the both real parts
$X_\sigma(\r), X_{\tau\sigma }(\r)$ are non-empty.
Then we give a formula connecting the number $s_{nor}$ with
some invariants of the action of the group $\Gamma $ on the lattice
$H^2(X(\cc);\z)$ with the intersection pairing.
We mention that it is not clear that one can express the number
$s_{nor}$ using the action of $\Gamma$ on the lattice
$H^2(X(\cc);\z)$. This formula is very important for the
topological classification of real Enriques surfaces (see
\cite{N4}). See \S 3 for details.

I am grateful to O. Gabber for assuring me
that the statement of Proposition 1.1 below should be true.
I am grateful to R. Sujatha for very useful
discussions, in particular, for pointing out me on the results of E. Witt
from \cite{W}.

A preliminary variant of this paper \cite{N3}
was written during my stay
in the University of Notre Dame (USA) at 1991---1992.
I am grateful to the University of Notre Dame for hospitality.

\head
\S 1. The proof of the Theorem 0.1.
\endhead

We recall that if $X$
is a topological space with an action of a group $G$  and $\Cal A$ is a
$G$-sheaf of groups on $X$, then the group $H^k(X;G,{\Cal A})$
of equivariant cohomology (or Galois cohomology) is defined.
See A. Grothendieck \cite{Gr1, Ch. 5}.
It is the right derived
functor $R^k\Gamma ^G$ to the  functor
${\Cal A}\mapsto \Gamma (X;{\Cal A})^G$
of $G$-invariant sections. This composition of functors
${\Cal A}\mapsto \Gamma (X;{\Cal A})$ and
$M\mapsto M^G$ for a $G$-module $M$
defines two spectral functors which tend to this cohomology:
$$
I_2^{p,q}=H^p(X/G;{\Cal H}^q(G;{\Cal A}))\Longrightarrow
H^{p+q}(X;G,{\Cal A})
$$
where ${\Cal H}^q(G;{\Cal A})=R^qf_\ast ^G{\Cal A}$ is the sheaf
corresponding to the presheaf  on $X/G$:
$$
U \longmapsto H^q(\pi ^{-1}(U);G, {\Cal A}).
\tag1--1
$$
Here $\pi :X\to X/G$ is the quotient map.
And
$$
II_2^{p,q}=H^p(G;H^q(X;{\Cal A})) \Longrightarrow H^{p+q}(X;G,{\Cal A}).
 \tag1--2
$$

The following statement is fundamental for us. This is analogous to the
well-known connection between \'etale and ordinary cohomology for
a complex algebraic manifold $Z$ and a finite abelian group $B$ (see
\cite{Mi, Ch. III, \S 3}):
The morphism of the ordinary site (Euclidean topology) to
the \'etale site gives rise an isomorphism
$$
\Het^k(Z;B)\cong H^k(Z(\cc );B).
\tag1--3
$$

\proclaim{Proposition 1.1} Let $X/\r $ be a real algebraic manifold and
$G=Gal(\cc/\r )$. Let ${\Cal A}$ be a sheaf in \'etale topology of $X$
such that this sheaf is a constant $G$-sheaf
$A$ on $X\otimes \cc$, where $A$ is a finite
abelian group with an action of the group $G$.

Then there exists the canonical isomorphism
$\Het^k(X;A)\cong H^k(X(\cc);G,A)$
together with the canonical isomorphism of
the Hochschild-Serre spectral sequence and
the Grothendieck
spectral sequence $II$, which is defined by the canonical isomorphism
$$
E_2^{p,q}=H^p(G;\Het^q(X\otimes \cc;A)) \cong
II_2^{p,q}=H^p(G;H^q(X(\cc);A))
$$
induced by the canonical isomorphism (1--3).
\endproclaim

\demo{Proof} We follow to the proof of the isomorphism (1--3)
for complex algebraic manifolds. See \cite{Mi,Ch. III, \S 3}, for
example.

Let $X(\cc)_{cx}$ be a small site $X(\cc )^{an}_E$
of morphisms of complex analytic spaces over $X(\cc )^{an}$,
which are local isomorphisms. (We use notation of
\cite{Mi}.) An open subset $U\subset X(\cc )$ is a local isomorphism. It
follows that we have the morphism of sites
$X(\cc)_{cx}\to X(\cc )$. Every covering of $X(\cc )$ in the site
$X(\cc)_{cx}$ has a refinement covering of $X(\cc )$ in the site
$X(\cc)$. It follows that
the morphism $X(\cc)_{cx}\to X(\cc)$ gives an isomorphism of cohomology
$$
H^i(X(\cc);A)\cong H^i(X(\cc)_{cx};A).
\tag1--4
$$
Like above, we can define
$G$-equivariant cohomology
$H^k(X(\cc)_{cx};G,{\Cal F})$ for the site $X(\cc)_{cx}$ as the
right derived functor to the composition of functors
${\Cal F}\to \Gamma (X(\cc )_{cx};{\Cal F})^G$ for a
$G$-sheaf ${\Cal F}$ on the site $X(\cc)_{cx}$.
This equivariant cohomology also has
a spectral sequence $II(X(\cc)_{cx})$ with the beginning
$$
II_2^{p,q}(X(\cc)_{cx})=
H^p(G;H^q(X(\cc)_{cx};{\Cal F})) \Longrightarrow
H^{p+q}(X(\cc)_{cx};G,{\Cal F}).
$$
One can calculate equivariant cohomology using invariant
coverings (see \cite{Gr1, Ch.5}).
Any $G$-invariant covering
of $X(\cc )$ in site $X(\cc)_{cx}$ also contains
a refinement $G$-invariant covering of
$X(\cc)$ in the site $X(\cc)$. It follows that
we also have the canonical isomorphism of equivariant cohomology:
$$
H^i(X(\cc)_{cx};G,A)\cong H^i(X(\cc);G,A),
\tag1--5
$$
together with the isomorphism of the
corresponding spectral sequences $II$ of these cohomology
$$
II_2^{p,q}(X(\cc)_{cx})=H^p(G;H^q(X(\cc)_{cx};A))\cong
II_2^{p,q}(X(\cc))=H^p(G;H^q(X(\cc);A)).
\tag1--6
$$
defined by the isomorphism (1--4).

Further, we use standard results about \'etale cohomology (see
\cite{Mi}). Since $X\otimes \cc \to X$ is an \'etale
covering with Galois group $G$, it follows that a sheaf
$F$ on the site $X_{et}$ corresponds to a
$G$-sheaf $F$ on the site $(X\otimes \cc)_{et}$.
\'Etale cohomology is a right derived
functor to the composition of
functors
$$
F \to (G-mod\ \ \Gamma (X\otimes \cc ; F))\ \ \text{and} \ \
M\to M^G.
$$
Here $M$ is a $G$-module.
The Hochschild-Serre spectral sequence
$$
E_2^{p,q}=H^p(G;\Het^q(X\otimes \cc; F))\Longrightarrow
\Het^{p+q}(X; F)
$$
corresponds to the composition of these functors.

Every \'etale  morphism $Y\to X\otimes \cc$ gives
a morphism $Y(\cc) \to X(\cc)$ in the site $X(\cc )_{cx}$.
This defines the morphism of sites
$X(\cc)_{cx} \to X_{et}$. From the remarks above, this morphism
defines the homomorphism of cohomology
$$
\Het^k(X;A)\to H^k(X(\cc)_{cx};G,A)=H^k(X(\cc);G,A)
\tag1--7
$$
together with the homomorphism of spectral sequences
$$
E_r^{p,q}\to II_r^{p,q}
\tag1--8
$$
defined by the homomorphism
$$
\split
E_2^{p,q}=H^p(G; \Het^q(X\otimes \cc;A))\to
II_2^{p,q}&=H^p(G;H^q(X(\cc)_{cx};A))\\
&=H^p(G;H^q(X(\cc);A)).
\endsplit
\tag1--9
$$
The last homomorphism is defined by the homomorphism
$$
\Het^q(X\otimes \cc; A)\to H^q(X(\cc )_{cx};A)=H^q(X(\cc);A),
\tag1--10
$$
which is the isomorphism (1--3).
It follows that (1--9), (1--8) and (1--7) are isomorphisms too.
This finishes the proof.

We mention that D. A. Cox \cite{C} had shown that \'etale homotopy
type of a real algebraic manifold is defined by Euclidean topology.

\enddemo

We recall that for a compact manifold $M$ and a constant
sheaf of modules, one can calculate sheaf
cohomology using simplicial triangulation
of the manifold $M$. One
can prove this using the canonical isomorphisms between
sheaf, \v Chech, Alexander--Spanier, singular and
simplicial triangulation cohomology for compact manifolds.
See \cite{Gr1} and \cite{Sp}, for example.
Similarly, for a compact manifold $M$ with an action of the group $G$ and
an Abelian $G$-group $A$, one can prove that equivariant
sheaf, \v Chech, Alexander--Spanier, singular
and simplicial triangulation cohomology are isomorphic.
Thus, for this case, we can calculate
equivariant cohomology using the following elementary procedure (this
definition is used in the book \cite{Bro}, for example):

We consider  some $G$-equivariant simplicial triangulation $K$ of $M$.
Thus, $K$ is a
$G$-equivariant simplicial complex. We may suppose that
$K/G$ is a simplicial triangulation of $M/G$ and the fixed part
$K^G$ is a simplicial triangulation of $M^G$. We consider the
corresponding chain complex $C_n=C_n(K;\z)$ and the corresponding cochain
complex $C^n(K;A)$, where $C^n(K;A)=\text{Hom}(C_n,A)$.
We organize a cochain complex for calculation of group cohomology of
$C^n=C^n(K,A)$.
For the group $G=\{ 1,g \} $ of order two this is a complex
$$
0\to C^n @> 1-g >>C^n @> 1+g >>C^n @> 1-g >>
C^n @>1+g>> C^n @>1-g>> \cdots
\tag1--11
$$
Thus, we get a double complex
$$
\CD
\vdots @. \vdots @. \vdots   @.  \vdots
@.   \vdots    @. \vdots  @.   \cdots \\
@.     @Ad^4AA              @Ad^4AA             @Ad^4AA
        @Ad^4AA            @Ad^4AA \\
0 @>>> C^4  @> 1-g  >> C^4  @> 1+g >> C^4 @>
1-g >>  C^4  @>1+g>>  C^4 @>1-g>>\cdots \\
@.     @Ad^3AA              @Ad^3AA             @Ad^3AA
        @Ad^3AA            @Ad^3AA \\
0 @>>> C^3  @> 1-g  >> C^3  @> 1+g >> C^3 @>
1-g >>  C^3  @>1+g>>  C^3 @>1-g>>\cdots \\
@.     @Ad^2AA              @Ad^2AA             @Ad^2AA
        @Ad^2AA            @Ad^2AA \\
0 @>>> C^2  @> 1-g  >> C^2  @> 1+g >> C^2 @>
1-g >>  C^2  @>1+g>>  C^2 @>1-g>>\cdots \\
@.     @Ad^1AA              @Ad^1AA             @Ad^1AA
        @Ad^1AA            @Ad^1AA \\
0 @>>> C^1  @> 1-g  >> C^1  @> 1+g >> C^1 @>
1-g >>  C^1  @>1+g>>  C^1 @>1-g>>\cdots \\
@.     @Ad^0AA              @Ad^0AA             @Ad^0AA
        @Ad^0AA            @Ad^0AA \\
0 @>>> C^0  @> 1-g  >> C^0  @> 1+g >> C^0 @>
1-g >>  C^0  @>1+g>>  C^0 @>1-g>>\cdots \\
@.     @AAA                 @AAA                @AAA
        @AAA            @AAA \\
@.      0       @.        0         @.        0
@.     0     @.         0
\endCD
\tag1--12
$$
Equivariant cohomology $H^n(M;G,A)$ are the homology of this
double complex. The spectral sequence $I$
corresponds to the filtration from
below on cohomology of this double complex, and the spectral sequence $II$
corresponds to the filtration from left of this double complex. From this
complex, it is clear (see Godement \cite{Go, Chapter 1} , for example) that
$$
II_1^{p,q}(M;G,A)=H^q(M;A)\Longrightarrow H^{p+q}(M;G, A),
\tag1--13
$$
and
$$
II_2^{p,q}(M;G,A)=H^p(G;H^q(M;A)) \Longrightarrow
H^{p+q}(M;G,A).
\tag1--14
$$
Here, for a $G$-module $R$, we have
$$
H^0(G;R)=R^G,
\tag1--15
$$
$$
 H^p(G;R)=R^g/(1+g)R \ \text{for}\ p \  \text{even and } p>0;
\tag1--16
$$
and
$$
H^p(G;R)=R^{(-g)}/(1-g)R \ \text{for}\  p \  \text{odd}.
\tag1--17
$$
If $R$ is a vector space over the field $\f _2$, the formulae
(1--16) and (1--17) give the same.

Further, we calculate the
equivariant cohomology for $A=\f _2$. For this case,
$1-g=1+g$, and the double complex (1--12) is periodic.

The spectral sequence $I$ corresponds to the filtration of this complex
from below.
Then
$$
I_1^{p,0}(M;G,\f_2) = (C^p)^g=C^p(M/G;\f _2),
\tag1--18
$$
and
$$
I_2^{p,0}(M;G,\f _2)=(C^p)^g/d^{p-1}((C^{p-1})^g)=H^p(M/G;\f _2).
\tag1--19
$$
For $q>0$,
$$
\split
I_1^{p,q}(M;G,\f_2) = (C^p)^g/(1+g)(C^p) & =
C^p(M/G;\f_2)/C^p(M/G,M^G;\f_2)\\
& =C^p(M^G;\f_2).
\endsplit
\tag1--20
$$
And we get
$$
I_2^{p,q}(M;G,\f_2) =H^p(M^G;\f _2)\ \  \text{if} \ \ q>0.
\tag1--21
$$
In particular, for a real projective algebraic manifold $X$,
the set $X(\cc)$ of complex points with the group $G=\{1,g\}$
generated by antiholomorphic involution $g$ acting on
$X(\cc)$, and the set
$X(\r)=X(\cc )^G$ of real points, we have:
$$
II_2^{p,q}(X(\cc);G,A)=H^p(G;H^q(X(\cc);A)) \Longrightarrow
H^{p+q}(X(\cc);G,A);
\tag1--22
$$
$$
I_2^{p,0}(X(\cc);G,\f _2)=H^p(X(\cc)/G;\f _2);
\tag1--23
$$
and
$$
I_2^{p,q}(X(\cc);G,\f_2) =H^p(X(\r);\f _2)\ \  \text{if} \ \ q>0.
\tag1--24
$$
Further, we will consider this case of real projective
algebraic manifold $X$, but, actually, all results
valid for a compact manifold $M$ with an involution.

For $X(\r )$, the double complex (1--12) has
$1-g=1+g=0$. Thus, we evidently get
(see \cite{Gr1, Corollaire 5.4.1})

\proclaim{Proposition 1.2} We have a canonical isomorphism
$$
H^k(X(\r );G, \f _2)=
\bigoplus_{i=0}^k {H^i(X(\r );\f _2)}.
$$
Besides, the spectral sequence $I(X(\r );G,\f _2)$ degenerates
from $I_2$, and for all $p,q$ we have
$I^{p,q}_2(X(\r);G,\f_2)=H^p(X(\r );\f _2)$.
Thus, all differentials
$$
d_r^{p,q}(X(\r );G,\f _2)$$
of the spectral sequence $I$ vanish for $r\ge 2$.
\endproclaim

The same is true for a topological space with the trivial action of the
group $G$.

Now we have (see \cite{Kr})

\proclaim{Proposition 1.3} For $k>2\dim~X$,
$$
H^k(X(\cc );G, \f _2)=
H^k(X(\r );G, \f _2)=\bigoplus_{i=0}^k H^i(X(\r );\f _2).
$$
\endproclaim

\demo{Proof} Like above, using the double complex, we can define
equivariant cohomology of a pair. We evidently have an exact sequence;
$$
\hdots \to H^{i}(X(\cc), X(\r );G,\f _2)\to
H^i(X(\cc);G,\f _2)\to H^i(X(\r );G,\f _2) \to
$$
$$
H^{i+1}(X(\cc), X(\r );G,\f _2)\to \hdots .
$$
The group $G$ acts without fixed points on the pair $(X(\cc ),X(\r ))$. --
For the corresponding double complex all differentials $1+g,1-g$ give an
exact sequence. It follows that for $q>0$ the spectral sequence
$I_1^{p,q}=0$ for this double complex and
$$
H^i(X(\cc), X(\r );G,\f _2)=H^i(X(\cc)/G, X(\r );G,\f _2).
$$
Thus,
$H^i(X(\cc), X(\r );G,\f _2)=0$ for $i>2\dim~X$.
It follows the statement.
\enddemo

 From the calculation above of the beginning of the spectral sequence $I$ and
Proposition 1.2, we get

\proclaim{Proposition 1.4} The embedding $\rho:X(\r )\subset X(\cc )$ gives
the homomorphism of the spectral sequences
$\rho^\ast :I(X(\cc );G,\f_2)\to I(X(\r );G,\f _2)$. For
$r=2$ it is defined by
$$
I_2^{p,0}(X(\cc );G,\f_2)=H^p(X(\cc )/G;\f_2)@>i ^\ast >>H^p(X(\r );\f_2)=
I_2^{p,0}(X(\r );G,\f_2),
$$
where $i:X(\r )\subset X(\cc )/G$ is the embedding,
and by the identical isomorphism
$$
I_2^{p,q}(X(\cc );G,\f_2)=H^p(X(\r );\f_2)=I_2^{p,q}(X(\r );G,\f_2)
\ \ \text{for}\ \ q>0.
$$
In particular, since the spectral sequence $I(X(\r );G,\f_2)$ degenerates,
the differential
$d_r^{p,q}$ of
$I(X(\cc );G,\f_2)$ vanish if $r\ge 2$ and $q\not=r-1$, and
$\rho ^\ast d_r^{p,r-1}=0$ for $r\ge2$.
\endproclaim

Now we can reformulate the map (0--1)
$$
_2Br^\prime (X)\to (\z /2)^s.
$$
Let us choose a point $P_i$, $i=1,2,...,s$, for
every connected component of $X(\r )$.
Here $i=1,2,\cdots , s$ numerates connected components of
$X(\r )$.  From the Kummer exact sequence (0--3), from the isomorphism
$Br^\prime (\{ P_i\} )=H^2_{et}(\{P_i\} ;\f _2)=\f _2$
and Proposition 1.1 (also see Proposition 1.2),
the image of the homomorphism  (0--1) is the same as the image of the
composition of the homomorphisms
$$
H^2(X(\cc );G,\f _2)\to H^2(X(\r);G,\f_2)\to
\bigoplus _{i=1}^s {H^2(\{ P_i\} ;G,\f_2)}.
\tag1--25
$$
defined by the inclusions $\{ P_i\} \subset X(\r )\subset X(\cc )$.
By Proposition 1.2, the second homomorphism is epimorphic and has
the kernel
$$
H^1(X(\r );\f _2)\oplus H^2(X(\r );\f _2).
$$
-- It is clear from the definition of equivariant cohomology using the
double complex (1--12).
We remark that we then have a canonical identification
$$
\split
H^0(X(\r );\f _2)&=H^2(X(\r );G,\f _2)/(H^1(X(\r );\f _2)\oplus
H^2(X(\r );\f _2)\\
&=I^{0,2}_\infty (X(\r );G,\f _2)
\endsplit
$$
for the first spectral sequence $I$ of the equivariant cohomology of
$X(\r )$.  It follows, that the image of the map
(1--25) is the image of the canonical homomorphism
$$
H^2(X(\cc);G,\f_2)\to H^2(X(\r);G,\f_2) \to
I_\infty^{0,2}(X(\r);G;\f _2)
=H^0(X(\r);\f_2).
$$
This homomorphism preserves the filtration $I$ on $H^2(X(\cc );G,\f_2)$.
Thus, this image is the same as the image of the canonical homomorphism
$$
\rho^\ast : I_\infty^{0,2}(X(\cc );G,\f _2)\to I_\infty^{0,2}(X(\r );G,\f _2)=
H^0(X(\r );\f_2),
$$
where $\rho :X(\r )\subset X(\cc )$ is the embedding.

Let us look on the part of the spectral sequence $I_2^{p,q}(X(\cc );G,\f _2)$
which  takes part in calculation of the $I_\infty ^{0,2}(X(\cc );G,\f _2)$.
It is the left-below corner
$$
\alignat 3
&H^3(X(\cc )/G;\f _2)\\
&H^2(X(\cc )/G;\f _2) &\qquad  &H^2(X(\r );\f _2)\\
&H^1(X(\cc )/G;\f _2) & \qquad &H^1(X(\r );\f _2) &\qquad & H^1(X(\r );\f _2)\\
&H^0(X(\cc )/G;\f _2) & \qquad &H^0(X(\r );\f _2) &\qquad &H^0(X(\r );\f _2)\\
\endalignat
$$
where the differential hits from below to up and from right to left.
In particular, we have $I_\infty ^{0,2}(X(\cc );G,\f _2)\subset
I_2^{0,2}(X(\cc );G,\f _2)=H^0(X(\r);\f _2)$. From the double
complex (1--12), the canonical map
$$
\split
\rho ^\ast :
I_2^{0,2}(X(\cc );G,\f_2)&=H^0(X(\r );\f _2) \\
&\to
I_2^{0,2}(X(\r );G,\f _2)=I_\infty ^{0,2}(X(\r );G,\f _2)=H^0(X(\r );\f _2)
\endsplit
$$
is the identical isomorphism of $H^0(X(\r );\f _2)$.
Thus, we get

\proclaim{Proposition 1.5} The image of the map $\rho ^\ast $
(equivalently, of the map (0--1))
is equal to the image of the embedding
$$
I_\infty ^{0,2}(X(\cc );G,\f_2) \subset I_2 ^{0,2}(X(\cc );G,\f_2)=
H^0(X(\r );\f_2)= I_\infty^{0,2}(X(\r );G,\f_2).
$$
\endproclaim

 From Proposition 1.4, the differential
$$
d_2^{0,2}:H^0(X(\r);\f_2)\to H^2(X(\r);\f_2)
$$
is equal to zero. By Proposition 1.4, we also have
$i^\ast d_2^{1,1}=0$ and $\rho^\ast d_3^{0,2} =0$.
We then get the basic result

\proclaim{Theorem 1.6} Let $i:X(\r)\subset X(\cc )/G$ be the embedding and
\linebreak
$i^\ast :H^3(X(\cc )/G;\f _2)\to H^3(X(\r );\f_2)$ the corresponding
canonical homomorphism.

Then we have:
The differential $d_2^{1,1}$ is the composition
$$
d_2^{1,1}:H^1(X(\r);\f_2)\to \Ker i^\ast \subset H^3(X(\cc )/G;\f _2),
$$
and the image of the homomorphism $\rho ^\ast $ from Proposition 1.5
(equivalently, of the homomorphism (0--1) from Introduction) is  equal to the
kernel of the differential
$$
\split
d_3^{0,2}:H^0(X(\r );\f _2) &\to \Ker i^\ast /d_2^{1,1}(H^1(X(\r);\f_2))\\
&\subset
H^3(X(\cc )/G;\f _2)/d_2^{1,1}(H^1(X(\r);\f _2)).
\endsplit
$$
In particular, $\rho ^\ast $ is epimorphic if $\Ker i^\ast =0$ or
$H^3(X(\cc )/G;\f _2)=0$.
It is epimorphic too if the differential
$$
d_2^{1,1}:H^1(X(\r);\f _2) \to \Ker i^\ast
$$
is epimorphic.
\endproclaim

As a corollary, we get the Theorem 0.1 from Introduction.

We remark that for
surfaces $X$ the dimension $\dim~X(\r )=2$, and \linebreak
$H^3(X(\r );\f_2)=0.$ Thus,
$\Ker i^\ast =H^3(X(\cc )/G;\f_2)$.

\remark{Remark 1.7} In Introduction (after Corollary 0.2),
we showed real
Enriques surfaces  $Y$ such that
$$
H^3(Y(\cc )/G;\f_2)=\f _2.
$$
Theorem 1.6 shows that such a $Y$
gives an  example when the map (0--1) is not epimorphic, exactly if
simultaneously the homomorphism
$$
d_2^{1,1}:H^1(Y(\r );\f _2) \to H^3(Y(\cc)/G;\f_2)=\f_2
$$
is zero, and  the homomorphism
$$
d_3^{0,2}:H^0(Y(\r );\f _2) \to H^3(Y(\cc )/G;\f_2)=\f_2
$$
is not zero.

Of course, these homomorphisms are "the same" for Enriques surfaces which
belong to one connected component of the moduli space of real Enriques
surfaces.  Using Global Torelli Theorem \cite{P\u S-\u S}, epimorphisity of
Torelli map  \cite{Ku} for $K3$-surfaces and methods developed in
\cite{N1,N2}, in  principle, it is possible to enumerate all this connected
components using  some invariants.
Thus, the problem is to rewrite, using these invariants, the invariants of
differentials $d_2^{1,1},d_3^{0,2}$ above. We hope to do it later. Some
important results in this direction were obtained in \cite{N-S}.
\endremark

\remark{Remark 1.8} We can consider the homomorphism
$$
H^1_{et}(X;\f _2)\to (\z /2)^s
\tag1--26
$$
which is defined like the map (0--1) using the isomorphism
$H^1_{et}(\{ real\  point\} ;\f _2)=\f _2$.

Like above, we can interpret the right side of (1--26) as
the group \linebreak
$H^0(X(\r );\f _2)=I_{\infty }^{0,1}(X(\r );G,\f _2)$,
and the image of
the homomorphism (1--26) as the
$$
\Ker \{ d_2^{0,1}:H^0(X(\r ),\f _2)\to H^2(X(\cc )/G;\f _2) \}
$$
for the differential $d_2^{0,1}$ of the spectral sequence
$I(X(\cc );G,\f _2)$.

Suppose that $X$ is a smooth projective real curve and
$X(\r )\not=\emptyset$. Then the quotient space $X(\cc )/G$ is a connected
2-dimensional manifold with a not-trivial boundary $X(\r )$.
By Poincar\'e duality,
$H^2(X(\cc )/G;\f_2)=H_0(X(\cc)/G,X(\r);\f_2)=0$,
and the map (1--26) is epimorphic.
This gives a geometrical interpretation of the result of E.Witt from
\cite{W}.

We mention, that similarly to the homomorphisms
(0--1) and (1--26), one may define and
study the  general homomorphism  (1--27) below
$$
H_{et}^n(X;\f_2)\to (\z/2)^s
\tag1--27
$$
using the isomorphism
$H^n_{et}(\{ real\  point\} ;\f _2)=\f _2$ for $n\ge 0$.
\endremark

\head
\S 2. The Proof of the Theorems 0.3 --- 0.6.
\endhead

We prove the following

\proclaim{Theorem 2.1} Let $X/\r$ be a smooth real projective algebraic
surface such that
$X(\r )\not=\emptyset$
and $H^3(X(\cc )/G;\f_2)=0$.

Then the Serre-Hochschild spectral
sequence  for the $H_{et}^\ast (X;\f _2)$
and the spectral sequence $II(X(\cc);G,\f_2)$ for the
equivariant cohomology $H^\ast (X(\cc );G,\f _2)$ degenerate.
Besides (by Proposition 1.1),
$H^k_{et}(X;\f_2)\cong H^k(X(\cc);G,\f_2)$, and we have the formulae for
their dimensions:
$$
\dim~H^0(X(\cc );G,\f _2)=1;
$$
$$
\dim~H^1(X(\cc );G,\f _2)=\dim~H^1(X(\cc );\f _2)+1;
$$
$$
\dim~H^2(X(\cc );G,\f _2)=
\dim~H^2(X(\cc );\f _2)^G+\dim~H^1(X(\cc );\f _2)+1;
$$
$$
\split
\dim~H^3(X(\cc );G,\f _2)=
2\dim~H^2(X(\cc );\f _2)^G&-\dim~H^2(X(\cc );\f _2)\\
                          &+2\dim~H^1(X(\cc );\f _2)+1;
\endsplit
$$
$$
\split
\dim~H^k(X(\cc );G, \f _2)=
2\dim~H^2(X(\cc );\f _2)^G&-\dim~H^2(X(\cc );\f _2)\\
                          &+2\dim~H^1(X(\cc );\f _2)+2
\endsplit
$$
for $k\ge 4$.
\endproclaim

\demo{Proof} By Proposition 1.1, we should prove that the
spectral sequence \linebreak
$II(X(\cc );G,\f_2)$ degenerates. To prove this,
we use the following important result of V. A. Krasnov \cite{Kr}.

\proclaim{Proposition 2.2}
Let $X/\r $ be a real projective algebraic manifold.

Then the spectral sequence $II(X(\cc );G,\f_2)$ with the beginning
$$
II_2^{p,q}=H^p(G;H^q(X(\cc ) ,\f _2)) \Longrightarrow H^{p+q}(X(\cc );G,\f _2)
$$
degenerates iff
$$
\dim~H^\ast (X(\r );\f _2)=\dim~H^1(G;H^\ast (X(\cc );\f _2)).
$$
\endproclaim

\demo{Proof} Let $k>2\dim~X$. By Propositions 1.2 and 1.3,
$$
\dim~H^\ast (X(\r );\f _2)=\dim~H^k(X(\cc);G,\f _2).
$$
On the other hand,
$$
\bigoplus _{p+q=k} {II_2^{p,q}}=
\bigoplus _{p+q=k} {H^p(G;H^q(X(\cc );\f _2))} =
H^1(G;H^\ast (X(\cc );\f _2)).
$$
Thus, we have the inequality
$$
\dim~H^\ast (X(\r );\f _2)=
\dim~H^k(X(\cc);G,\f_2)
\le \dim~H^1(G;H^\ast (X(\cc );\f _2)).
\tag2--1
$$
Besides, (2--1) gives
an equality iff $II_2^{p,q}=II_\infty ^{p,q}$ for
$p+q=k>2\dim~X$. Thus, if the spectral sequence
$II(X(\cc );G,\f_2)$ degenerates,
we have the equality of the Proposition.

Now we assume that the equality of Proposition holds.
We have then proven that all
differentials  $d_r^{p,q}$ vanish for $r\ge 2$ and $p+q>2\dim~X$.
 From periodicity of the double complex (1--12), it follows that all
differentials $d_r^{p,q}$ vanish for $r\ge 2$.
\enddemo

We recall the Smith exact sequence for an action of a group $G=\{ 1,g\} $ of
order two (see \cite{Bre, Ch. III, \S 3}, for example).

We have the exact sequence for the chain complex
$C_n(K;\f_2)$ above with an action of $G$:
$$
0\to C_n(K;\f_2)^G\to C_n(K;\f_2) \to (1+g)(C_n(K;\f_2)) \to 0.
$$
Here we have the canonical identifications
$$
(1+g)(C_n(K;\f_2))=C_n(K/G;\f_2)/C_n(K^G;\f_2)=C_n(K/G,K^G;\f_2)
$$
and
$$
C_n(K;\f _2)^G=(1+g)C_n(K;\f _2)\oplus C_n(K^G;\f _2)=
C_n(K/G,K^G;\f _2)\oplus C_n(K^G;\f _2).
$$
Thus, we get the exact sequence
$$
0\to C_n(K/G,K^G;\f _2)\oplus C_n(K^G;\f _2)\to C_n(K; \f _2) \to
C_n(K/G,K^G;\f _2)\to 0.
$$
This gives the corresponding homological Smith exact sequence:
$$
\split
\hdots
& @> \delta _{n+1} >> H_n(K/G, K^G;\f _2)\oplus H_n(K^G;\f _2) @>i_n >>
H_n(K;\f _2) @>\rho _n >> H_n(K/G, K^G;\f _2) \\
& @>\delta _n >> H_{n-1}(K/G,K^G;\f _2)\oplus H_{n-1}(K^G;\f _2)
@> i_{n-1} >> \hdots .
\endsplit
$$
Thus, for a real algebraic variety, we have the homological
Smith exact sequence
$$
\hdots
 @> \delta _{n+1} >> H_n(X(\cc )/G,X(\r );\f _2)\oplus H_n(X(\r  );\f _2)
@>i_n >> H_n(X(\cc );\f _2) @>\rho _n >>
$$
$$
H_n(X(\cc )/G,X(\r  );\f _2) @>\delta _n >>
 H_{n-1}(X(\cc )/G,X(\r  );\f _2)\oplus H_{n-1}(X(\r  );\f _2)
\to \hdots
$$
We repeat some well-known standard facts connected with the Smith
exact sequence (see V.A.Rokhlin \cite{R}, for example).
 From this exact sequence, we get
$$
\split
\dim~H_{n-1}(X(\cc )/G,X(\r );\f _2) & +\dim~H_{n-1}(X(\r );\f _2)\\
& =\dim~\text{Im}\ \delta _n + \dim~\text{Im}\ i_{n-1}.
\endsplit
$$
Moreover,
$$
\split
\dim~\im \delta _n
&=\dim~H_n(X(\cc )/G,X(\r  );\f _2)-\dim~\im \rho _n \\
&=\dim~H_n(X(\cc )/G,X(\r  );\f _2)-\dim~H_n(X(\cc );\f _2)+\dim~\im i_n.
\endsplit
$$
Thus, we get
$$
\split
&\dim~H_{n-1}(X(\cc )/G,X(\r );\f _2)
+\dim~H_{n-1}(X(\r );\f _2)\\
&=\dim~H_n(X(\cc )/G,X(\r  );\f _2)-\dim~H_n(X(\cc );\f _2)\\
&+
\dim~\text{Im}\ i_n + \dim~\text{Im}\ i_{n-1}.
\endsplit
$$
Considering the sum by $n$, we then get
$$
\dim~H_\ast (X(\r );\f _2)
=2\sum _n {\dim~\text{Im}\ i_n} -\dim~H_\ast (X(\cc );\f _2).
\tag2--2
$$
 From the exact sequence,
$$
0\to H_\ast (X(\cc );\f _2)^G\to H_\ast (X(\cc );\f _2)\to
(1+g) H_\ast (X(\cc );\f _2)\to 0,
$$
we have
$$
\split
\dim~H^1(G;H_\ast (X(\cc );\f _2))
&=\dim~H_\ast (X(\cc );\f _2)^G -\dim~(1+g)H_\ast (X(\cc );\f _2)\\
&=2\dim~H_\ast (X(\cc ); \f _2)^G - \dim~H_\ast (X(\cc );\f _2).
\endsplit
\tag2--3
$$
It is clear that
$$
\text{Im}\ i_n\subset H_n(X(\cc );\f _2)^G.
\tag2--4
$$
Thus, from (2--2), (2--3) and (2--4), we have an inequality
$$
\dim~H_\ast (X(\r );\f _2) \le
\dim~H^1(G;H_\ast (X(\cc );\f _2)).
\tag2--5
$$
Of course, the inequalities (2--1) and (2--5) are equivalent.
Moreover, we see that the inequality (2--5) is an equality iff
for any $n$ we have for the Smith exact sequence
$$
\text{Im}\ i_n=H_n(X(\cc );\f _2)^G.
\tag2--6
$$
Thus, to prove that the spectral sequence $II$ degenerates,
we have to prove the equalities (2--6) for
$0\le n \le 4$ for a surface $X$
with the condition
$H^3(X(\cc )/G;\f _2)=0$.

It will be convenient for us using the following general statement which
follows from Smith exact sequence (compare with the proof of
\cite{H, Lemma 3.7}).

\proclaim{Proposition 2.3} For the Smith exact sequence,
$$
\split
\rho _n (H_n(X(\cc );\f_2)^G)
&=\text{Im}\
\{ H_{n+1}(X(\cc )/G;\f _2)\to H_{n+1}(X(\cc )/G,X(\r );\f _2)\\
&@>\delta _{n+1} >>
\text{Ker}\ \delta _n\subset H_n(X(\cc )/G,X(\r );\f _2)\}.
\endsplit
$$
In particular, $\text{Im}\ i_n=H_n(X(\cc );\f_2)^G$ iff the image
on the right is zero.
\endproclaim

\demo{Proof} We use the following properties of
the Smith exact sequence which follow from the definition above of
this sequence:
$$
i_n(\rho_n \oplus 0)=id+g;
\tag2--7
$$
and the homomorphism
$$
\split
H_n(X(\cc )/G;X(\r );\f _2)
&@>\delta _n>>
H_{n-1}(X(\cc )/G,X(\r );\f _2) \oplus H_{n-1}(X(\r );\f_2)\\
& @> \pi _{X(\r )} >> H_{n-1}(X(\r );\f _2)
\endsplit
\tag2--8
$$
is equal to the homomorphism
$\partial _n:H_n(X(\cc )/G,X(\r );\f _2)\to H_{n-1}(X(\r );\f _2)$
in the homological exact sequence of the pair
$(X(\cc )/G,X(\r ))$.

Now, let $x_n\in H_n(X(\cc );\f _2)^G$. Then, from (2--7), it is
equivalent to $i_n(\rho _n(x_n)\oplus 0)=0$. By Smith exact sequence,
it is equivalent to existence of an element
$y_{n+1}\in H_{n+1}(X(\cc )/G,X(\r );\f _2)$ such that
$\delta _{n+1} (y_{n+1})=\rho _n(x_n)\oplus 0$. By (2--8), it is equivalent
to
$$
\split
\rho_n (x_n)\in \text{Im}\
\{ H_{n+1}(X(\cc )/G;\f _2)
&\to H_{n+1}(X(\cc )/G,X(\r );\f _2)\\
&@>\delta _{n+1} >> H_n(X(\cc )/G,X(\r );\f _2)\}.
\endsplit
$$
Now, we should only remark that, by Smith exact sequence,
$\rho_n(x_n) \in \text{Ker}\  \delta _n$.

\enddemo

 From the Proposition 2.3, we have:

$\text{Im}\ i_4=H_4(X(\cc );\f _2)^G$
for any surface since $H_5(X(\cc )/G;\f _2)=0$.

$\text{Im}\ i_3=H_3(X(\cc );\f _2)^G$ since for our case $\Ker \delta _3=0$,
because $\Ker \delta _3\subset \Ker \partial _3=0$. Here
$\Ker \partial _3=0$, because $H_3(X(\cc )/G;\f _2)=0$ for our case.

$\text{Im}\ i_2=H_2(X(\cc );\f _2)^G$ since $H_3(X(\cc )/G;\f _2)=0$
in our case.

$\text{Im}\ i_1=H_1(X(\cc );\f _2)^G$ since
$\Ker \delta _1\subset \Ker \partial _1 =0$ in our case. Here
$\Ker \partial _1=0$ since
$H_1(X(\cc )/G;\f _2)=0$ in our case.

$\text{Im}\ i_0=H_0(X(\cc );\f _2)^G$ since $H_0(X(\cc ),X(\r );\f _2)=0$
because $X(\r )\not= \emptyset$.

Thus, we  proved that the spectral sequence $II$ degenerates.

Now let us prove the formulae of Theorems 1.2 and  0.3.
Since the spectral sequence  $II$ degenerates, we have
$$
\dim~H^k(X(\cc);G,\f_2)=\bigoplus _{p+q=k} {H^p(G;H^q(X(\cc );\f _2))}.
$$
To get formulae, by Poincar\'e duality,
we should only prove that $G$ is trivial
on \linebreak
$H^0(X(\cc );\f _2)$ and $H_1(X(\cc );\f _2)$. It is true for
$H^0(X(\cc );\f _2)$, since \linebreak
$H^0(X(\cc );\f _2)\cong \f _2$.
Since $H_1(X(\cc )/G;\f _2)=0$, the homomorphism $\partial _1$ is injective.
By (2--8), then the homomorphism $\delta _1$ for the Smith exact sequence is
injective too. Thus, the homomorphism $\rho _1$ is zero and
$H_1(X(\cc );\f _2)=\im i_1 \subset H_1(X(\cc );\f _2)^G$.
It follows the statement.
It finishes the proof of Theorems 2.1 and 0.3.
\enddemo

\demo{Proof of Theorem 0.4}
For $n\in \bold N$, the exact sequence of sheafs
$$
0 \to \z @>\times n >>\z \to \z /n \to 0
$$
gives the exact sequence of cohomology (universal coefficient sequence)
$$
\split
\cdots &\to  H^{k-1}(M;\z /n) \to H^k(M;\z )@>\times n>> H^k(M;\z ) \\
&\to H^k(M;\z /n)\to H^{k+1}(M;\z )@>\times n >> H^{k+1}(M;\z ) \to
\endsplit
\tag2--9
$$
For a compact manifold $M$,
the beginning of this sequence gives the
exact sequences
$$
0 \to H^0(M;\z )@>\times n>> H^0(M;\z ) \to H^0(M;\z /n)\to 0,
\tag2--10
$$
and
$$
0 \to H^1(M;\z )@>\times n>> H^1(M;\z )
\to H^1(M;\z /n)\to \hdots \ .
\tag2--11
$$
In particular, $H^1(M;\z)$ has no torsion.
As we had mentioned in Introduction, for a smooth surface $X$, the
quotient $X(\cc)/G$ is a smooth manifold. The
group $G$ preserves the canonical orientation of $X(\cc)$. It follows that
$X(\cc)/G$ is a smooth oriented manifold.
Since $H^3(X(\cc )/G;\f _2)=0$, by Poincar\'e duality \linebreak
$H^1(X(\cc )/G;\f _2)=0$. Thus, by (2--11), we then get that
$H^1(X(\cc )/G;\z )=0$. It follows
$H^1(X(\cc )/G;\cc )=0$. Since $2$ is invertible in $\cc $, we then
have $H^1(X(\cc );\cc )^G=\pi ^\ast H^1(X(\cc )/G;\cc )=0$ where
$\pi :X(\cc)\to X(\cc)/G$ is the quotient map
(see \cite{Bre, Ch. III, Theorem 2.4}).

For the  Hodge decomposition
$H^1(X(\cc);\cc )=H^{1,0}(X(\cc ))+H^{0,1}(X(\cc ))$,
the antiholomorphic
involution $g$ on  $X(\cc )$ evidently maps
$H^{1,0}(X(\cc ))$ to
$H^{0,1}(X(\cc))$.
It follows that $0=\dim
H^1(X(\cc );\cc )^G=(1/2)\dim~H^1(X(\cc );\cc )$.
Thus, applying (2--11) to
$X(\cc )$, we get
$$
H^1(X(\cc );\z )=0.
\tag2--12
$$
Thus, $X$ is a regular surface: the irregularity
$q(X)=\dim~H^{1,0}(X(\cc))=0$.
It follows that the characteristic class map gives an embedding
$$
\Pic (X\otimes \cc )\subset  H^2(X(\cc );\z ),
\tag2--13
$$
and the image of this map is defined by the condition
$$
\{ x\in \Pic (X\otimes \cc ) \ \mid \  x\cdot H^{2,0}(X(\cc))=0\} =
\gamma ^{-1}(H^{1,1}(X(\cc ))),
\tag2--14
$$
where
$$
\gamma :H^2(X(\cc );\z )\to H^2(X(\cc );\cc )
$$
is a coefficient map.
See \cite{G-H}, for example.
 From the definition of the characteristic class
map, we have
$$
g(\gamma (x))=-\gamma (g(x))\ \ \text{for}\ \ x\in \Pic (X\otimes \cc).
$$
We remark that $\Pic (X\otimes \cc )$ contains
the all torsion of $H^2(X(\cc );\z)$ by (2--14).

Since $H^1(X(\cc );\z )=0$ and the sequence
$$
0\to H^4(X(\cc );\z ) @>\times 2>> H^4(X(\cc );\z ) \to
H^4(X(\cc );\f _2) \to 0
$$
is exact,
from (2--9),
we get the exact sequence
$$
\split
0 &\to  H^1(X(\cc );\f _2) \to H^2(X(\cc );\z ) @>\times 2>> H^2(X(\cc );\z )
\\
&\to H^2(X(\cc );\f _2)\to H^3(X(\cc );\z ) @>\times 2>> H^3(X(\cc );\z ) \to
H^3(X(\cc );\f _2) \to 0.
\endsplit
\tag2--15
$$
Besides, by Poincar\'e duality,
$\dim~H^1(X(\cc );\f _2)=\dim~H^3(X(\cc );\f _2)$.
It follows that
$$
\dim~H^2(X(\cc );\z)/2H^2(X(\cc );\z )=b_2+\dim~H^1(X(\cc );\f _2),
\tag2--16
$$
and
$$
\dim~H^2(X(\cc );\f _2)=b_2+2\dim~H^1(X(\cc );\f _2)
\tag2--17
$$
where the Betti number $b_2=\dim~H^2(X(\cc );\cc )=
\dim~H^2(X(\cc );\z ) \otimes \cc $.

Let $\Pic (X\otimes \cc )= T\oplus \z ^{\rho(X\otimes \cc )}$
where $T$ is the
torsion of $\Pic (X\otimes \cc )$ and  $\z ^{\rho (X\otimes \cc )}=\Pic
(X\otimes \cc )/T$.
Since $X(\r )\not=\emptyset$, we have
$$
\Pic X= \Pic (X\otimes \cc )^G
$$
(this is well-known, see \cite{Ma}).
Let
$$
\Pic X= \Pic (X\otimes \cc )^G=T^\prime \oplus \z^{\rho (X)}
$$
where $T^\prime$ is the torsion of $\Pic X$ and
$\z^{\rho (X)}=\Pic X/T^\prime$.
If for $a\in \Pic (X\otimes \cc )$ we have
$g(a)=a\mod T$, then $g(ma)=ma$ for some $m \in \bold N$ such that
$mT=0$.
It follows that
$$
(\Pic (X \otimes \cc )/T)^G\cong (\Pic X)/T^\prime =\z ^{\rho (X)}.
$$
Thus,
$$
\rho (X)=\rho_+(X\otimes \cc ),
\tag2--18
$$
where $(\Pic (X \otimes \cc )/T)^G\cong \z^{\rho _+(X\otimes \cc )}$.
We had
proven above that
$G$ is trivial on $H_1(X(\cc );\f _2)$. Then, it is trivial
on $H^1(X(\cc );\f _2)=H_1(X(\cc );\f _2)^\ast$.
 From (2--15) and the remarks above, we  then get that the
group $G$ is trivial on
$$
\Ker \{ \Pic (X \otimes \cc ) @> \times 2 >> \Pic (X \otimes \cc )\} =
H^1(X(\cc );\f _2).
$$
Thus,
$$
\Ker \{ \Pic X  @> \times 2 >> \Pic X \} =
H^1(X(\cc );\f _2).
$$
As a result, we get that
$$
\dim~\Pic X/2\Pic X=\rho _+ (X\otimes \cc )+\dim~H^1(X(\cc );\f _2).
\tag2--19
$$
Here,
by the remarks above about the map (2--13),
$$
\rho_+ (X\otimes \cc )\le h_-^{1,1}(X(\cc )).
$$
 From Proposition 1.3 and formulae of Theorem 2.1, we have
$$
\split
\dim~H^\ast (X(\r );\f _2)& = \dim~H^5(X(\cc );G,\f_2)\\
& =2\dim~H^2(X(\cc );G,\f_2)-\dim~H^2(X(\cc );\f_2).
\endsplit
$$
Thus,
$$
\dim~H^2(X(\cc );G,\f_2)=(1/2)\dim~H^\ast (X(\r );\f _2)+
(1/2)\dim~H^2(X(\cc );\f_2).
\tag2--20
$$
By Proposition 1.1, $ H^2(X(\cc );G,\f_2)=
H^2_{et}(X;\f_2)$. By the exact sequence (0--3), we get
$\dim~_2Br^\prime (X)=\dim~H^2_{et}(X;\f_2)-\dim~\Pic X/2\Pic X$.
Thus, from (2--17), (2--19) and (2--20), we get
$$
\dim~_2Br^\prime (X)=
(1/2)\dim~ H^\ast (X(\r );\f _2)+b_2/2-\rho_+(X\otimes \cc ).
\tag2--21
$$
Let $(b_2)_+=\dim~H^2(X(\cc );\cc )^G$ and $(b_2)_-=b_2-(b_2)_+$.
 From the Lefschetz fixed point formula for the involution $g$
(see \cite{Sp}, for example), we get
$$
\chi(X(\r ))=2+2(b_2)_+ - b_2 = 2 + b_2 - 2(b_2)_-.
\tag2--22
$$
Thus, from (2--21) and (2--22), we get
$$
\dim~_2Br^\prime (X)=(1/2)\dim~ H^\ast (X(\r );\f_2)+(1/2)\chi(X(\r))-1
+(b_2)_- -\rho_+(X\otimes \cc ).
\tag2--23
$$
For the Hodge decomposition
$$
H^2(X(\cc );\cc)=H^{2,0}(X(\cc )) + H^{0,2}(X(\cc )) + H^{1,1}(X(\cc )),
$$
the antiholomorphic involution $g$ sends
$H^{2,0}(X(\cc ))\to H^{0,2}(X(\cc ))$ and \linebreak
$H^{1,1}(X(\cc ))\to H^{1,1}(X(\cc))$.
It follows that
$$
(b_2)_-=\dim~H^{2,0}(X(\cc ))+\dim~H^{1,1}_-(X(\cc )).
$$
For  a connected compact surface $F$ we have
$\dim~H^\ast (F;\f _2)+\chi (F)=4$.
Thus, from (2--23), we get the formula of Theorem 0.4.

\enddemo

\vskip 5pt

\demo{Proof of Theorem 0.6}  Let $Y$ be a real Enriques surface. In
\cite{N-S}, the inequality (0--7), i.e.
$b(Y)\ge 2s-2$, was proved. The proof was similar to the proof above
of the Theorem 0--4 and used Lefschetz fixed-point formula and the
inequality (2--1). From the proof, it follows that the
equality $b(Y)=2s-2$ holds iff
the inequality (2--1) is an equality.
By Proposition 2.2 (of V.A.Krasnov),
it then follows that the spectral sequence $II$ degenerates iff
$b(Y)=2s-2$. By Proposition 1.1, the Hochschild--Serre spectral
sequence degenerates iff $b(Y)=2s-2$. It follows the statement (i)
of Theorem 0.6. From the statement (i), the definition of the
invariant $\epsilon(Y)$, and from (0--6), (0--7),
the statements (ii) and (iii) of Theorem 0.6 follow.
\enddemo

\head
\S 3. Applications to topology of real Enriques surfaces
\endhead

We use notation on real Enriques surfaces of Introduction.
Thus, for a real Enriques surface $Y$, we denote by $\theta$
the antiholomorphic involution of $Y$, by $X$ the universal covering
$K3$-surface, by $\tau$ the holomorphic involution of the
$2$-sheeted universal covering $\pi :X(\cc)\to Y(\cc)$, and by
$\sigma $, $\tau\sigma$ two liftings of $\theta$ on $X(\cc)$. We
suppose that the automorphism group
$$
\Gamma = \{ \text{id},\ \tau ,\sigma ,\tau\sigma \}
$$
on $X(\cc)$
is isomorphic to $(\z/2)^2$. In particular, it is true if
$Y(\r)\not=\emptyset$ (see \cite{N-S}).

First, we discuss the following problem.
We have
$\dim~H^2(Y(\cc);\f_2)=12$. A subgroup
$H^2(Y(\cc);\z)\otimes \f_2\subset H^2(Y(\cc);\f_2)$ has
$\dim~H^2(Y(\cc);\z)\otimes \f_2=11$. Thus, we can introduce the
invariant
$$
\beta (Y)=\dim~H^2(Y(\cc);\f_2)^\theta -\dim
(H^2(Y(\cc);\z)\otimes \f_2)^\theta.
$$
This invariant is very important for real Enriques surfaces,
and first, we want to calculate $\beta(Y)$ in some cases.
Evidently, $\beta(Y)=0$ or $1$.
For the invariant
$$
b(Y)=
\dim~H^2(Y(\cc);\f_2)^\theta -
\dim~ (\Pic~Y\otimes \cc)^\theta /2(\Pic~Y\otimes \cc)^\theta + 1,
$$
(see Introduction), we have
$$
b(Y)=b^\prime (Y)+\beta(Y),
\tag3--1
$$
where we denote
$$
b^\prime (Y)=
\dim~(H^2(Y(\cc);\z)\otimes \f_2)^\theta -
\dim~ (\Pic~Y\otimes \cc)^\theta /2(\Pic~Y\otimes \cc)^\theta + 1.
$$
In \cite{N-S, Theorem 3.4.7}, there was obtained a formula for
$b^\prime(Y)$:
$$
b^\prime (Y)=r(\theta)-a(\theta)+\max\{1-\alpha(\sigma),\
(\delta_{\sigma L^{\tau ,\sigma}}+
\delta_{\sigma L^\tau_\sigma})/2\}.
\tag3--2
$$
Here $r(\theta), a(\theta),\alpha (\sigma),
\delta_{\sigma L^{\tau ,\sigma}},
\delta_{\sigma L^\tau_\sigma}$
are some invariants of the action of $\Gamma $ on
the lattice $L$ which is the lattice
$H^2(X(\cc);\z)$ with the intersection pairing. The invariants
$r(\sigma ), a(\sigma )$ are some non-negative integers, and
$$
r(\sigma )\equiv a(\sigma )\mod 2.
\tag3--3
$$
The invariants $\alpha (\sigma ) ,
\delta_{\sigma L^{\tau ,\sigma}},
\delta_{\sigma L^\tau_\sigma}$ are equal to $0$ or $1$,
and
$$
\delta_{\sigma L^{\tau , \sigma}}=\delta_{\sigma L^\tau_\sigma }.
\tag3--4
$$
The precise definition of these invariants is very long,
and we refer to \cite{N-S} for their definition. Actually,
these invariants are some specialization to Enriques
surfaces of general invariants from \cite{N1, N2a}.

Besides, in \cite{N-S}, it was proved that
the invariant $\beta (Y)=0$ if
$$
\max\{1-\alpha(\sigma),\
(\delta_{\sigma L^{\tau ,\sigma}}+
\delta_{\sigma L^\tau_\sigma})/2\}=0
$$
or, equivalently,
$\alpha (\sigma )=1$ and
$\delta_{\sigma L^{\tau ,\sigma}}=
\delta_{\sigma L^\tau_\sigma}=0$.

We want to prove here the following result which also gives
another prove of
the statement about $\beta(Y)$ above, but only in the case if
both real parts $X_\sigma(\r)$ and $X_{\tau\sigma}(\r)$ are non-empty.

\proclaim{Theorem 3.1} Let $Y$ be a real Enriques surface and
both $X_\sigma(\r)$ and $X_{\tau\sigma}(\r)$ are non-empty.

Then:
$$
\beta(Y)=\max\{1-\alpha(\sigma),\
(\delta_{\sigma L^{\tau ,\sigma}}+
\delta_{\sigma L^\tau_\sigma})/2\},
$$
and
$$
b(Y)=r(\theta)-a(\theta)+
2\max\{1-\alpha(\sigma),\
(\delta_{\sigma L^{\tau,\sigma}}+
\delta_{\sigma L^\tau_\sigma})/2\},
$$
where
$r(\theta),\ a(\theta),\
\alpha(\sigma),\
\delta_{\sigma L^{\tau,\sigma}},\
\delta_{\sigma L^\tau_\sigma}$
are some invariants of the action of the group $\Gamma $ on
the lattice $H^2(X(\cc);\z)$ with the intersection pairing.
\endproclaim

\demo{Proof} By Theorem 0.6, $b(Y)=2s-2\equiv 0 \mod 2$. By
(3--1) --- (3--3), we then get
$$
\beta(Y)\equiv
\max\{1-\alpha(\sigma),\
(\delta_{\sigma L^{\tau ,\sigma}}+
\delta_{\sigma L^\tau_\sigma})/2\} \mod~2.
$$
The right side of this congruence is equal to $0$ or $1$ since
$\alpha(\sigma)=0$ or $1$, and
$\delta_{\sigma L^{\tau ,\sigma}}=
\delta_{\sigma L^\tau_\sigma } = 0$ or $1$. It follows the first
formula since $\beta(Y)$ is equal to $0$ or $1$ too.
 From the first formula and (3--1), (3--2), the second one follows.
\enddemo

We don't know if this statement valid when
$\max\{1-\alpha(\sigma),\
(\delta_{\sigma L^{\tau ,\sigma}}+
\delta_{\sigma L^\tau_\sigma})/2\}=1$ and
one of $X_\sigma(\r)$ or $X_{\tau\sigma}(\r)$ is empty.

In \cite{N-S, Theorems 3.5.1---3.5.3, formula (3-5-1)},
there was obtained a formula for
$b(Y)$ using the numbers $s_{or}$ and $s_{nor}$ of orientable and
non-orientable connected components
of $Y(\r)$ respectively. The numbers $s_{or}$ and $s_{nor}$
are connected with the numbers $s(\sigma )$ and $s(\tau\sigma )$ of
connected components of $X_\sigma (\r )$ and
$X_{\tau\sigma}(\r)$ respectively
by the formula $s(\sigma )+s(\tau\sigma )=
2s_{or}+s_{nor}$ (see \cite{N-S, Lemma 3.2.1}).
For the case $s(\sigma )>0$ and $s(\tau\sigma )>0$
(or when both sets $X_\sigma (\r)$ and $X_{\tau\sigma}(\r )$ are
non-empty) which is necessary for us , this formula for $b(Y)$
claims that
$$
\split
b(Y)&=2s_{or}+s_{nor}-2 +\min~\{\alpha(\sigma),\
    (\delta_{\sigma L^{\tau ,\sigma}}+
    \delta_{\sigma L^\tau_\sigma})/2\}\\
& + \dim~H(\sigma )_- - \dim~H(\sigma)_+^\perp \cap H(\sigma)_- + \beta(Y).
\endsplit
$$
Here $\dim~H(\sigma )_-$ and $\dim~H(\sigma )_+^\perp \cap H(\sigma)_-$ are
some other invariants ot the action of the group $\Gamma$ on the lattice
$L$ (see \cite{N-S}).
 From the formula for $\beta (Y)$ of Theorem 3.1, we get the formula
$$
\split
b(Y)&=2s_{or}+s_{nor}-2+
\min\{\alpha (\sigma ),\
(\delta_{\sigma L^{\tau ,\sigma}}+
\delta_{\sigma L^\tau_\sigma})/2\}\\
& +\max\{1-\alpha(\sigma ),\
    (\delta_{\sigma L^{\tau ,\sigma}}+
   \delta_{\sigma L^\tau_\sigma})/2\}+
   \dim~H(\sigma )_- - \dim~H(\sigma )_+^\perp \cap H(\sigma )_- \\
    &= 2s_{or} + s_{nor} - 1 +
\alpha (\sigma )(\delta_{\sigma L^{\tau ,\sigma}}+
\delta_{\sigma L^\tau_\sigma}-1)\\
&+ \dim~H(\sigma )_- - \dim~H(\sigma )_+^\perp \cap H(\sigma )_-,
\endsplit
\tag3--5
$$
if both $X_\sigma (\r )$ and $X_{\tau\sigma}(\r )$ are non-empty.

By the formula $b(Y)=2s-2=2s_{or}+2s_{nor}-2$ of Theorem 0.6, we get

\proclaim{Theorem 3.2} Let $Y$ be a real Enriques surface and
both $X_\sigma (\r)$ and $X_{\tau\sigma}(\r )$ are non-empty.

Then for the number $s_{nor}$ of non-orientable connected components of
$Y(\r)$ we have the formula:
$$
s_{nor}=1 +
\alpha(\sigma )(\delta_{\sigma L^{\tau ,\sigma}}+
\delta_{\sigma L^\tau_\sigma} - 1)+
\dim~H(\sigma )_- -\dim~H(\sigma )_+^\perp \cap H(\sigma )_- ,
\tag3--6
$$
where $\alpha(\sigma ),\ \delta_{\sigma L^{\tau ,\sigma}},\
\delta_{\sigma L^\tau_\sigma},\
\dim~H(\sigma )_- ,\
\dim~H(\sigma )_+^\perp \cap H(\sigma )_- $ are some invariants of
the action of the group $\Gamma$ on the lattice $H^2(X(\cc);\z)$ with
the intersection pairing.

\endproclaim

We mention that by Theorem 0.6, $b(Y)=2s-2$ if both
$X_\sigma(\r)$ and $X_{\tau\sigma}(\r)$ are non-empty. Thus,
by the formula for $b(Y)$ of Theorem 3.1, we also have the formula
for the number $s=s_{or}+s_{nor}$ of all connected components of
$Y(\r)$ if both $X_\sigma(\r)$ and $X_{\tau\sigma}(\r)$ are non-empty.

\proclaim{Theorem 3.3} Let $Y$ be a real Enriques surface and
both $X_\sigma (\r)$ and $X_{\tau\sigma}(\r )$ are non-empty.

Then for the number $s$ of all connected components of $Y(\r)$ we
have the formula
$$
s=1+(r(\theta)-a(\theta))/2+
\max\{1-\alpha(\sigma),\ (\delta_{\sigma L^{\tau,\sigma}}+
\delta_{\sigma L^\tau_\sigma})/2\},
$$
where $r(\theta), \ a(\theta),\ \alpha(\sigma),\
\delta_{\sigma L^{\tau,\sigma}},\
\delta_{\sigma L^\tau_\sigma}$
are some invariants of the action of the group $\Gamma$ on the lattice
$H^2(X(\cc);\z)$ with the intersection pairing.
\endproclaim

Of course, from the formulae for $s$ and $s_{nor}$ or
Theorems 3.2, 3.3, we get a formula for $s_{or}=s - s_{nor}$.
These Theorems 3.1---3.3 are very important for the topological
classification of real Enriques surfaces---describing of all possible
topological types of $Y(\r)$ for real Enriques surfaces $Y$.
See \cite{N4} for these applications.

\enddocument

\end